\documentclass[11pt]{article}
\usepackage[a4paper,margin=27mm]{geometry}
\usepackage[T1]{fontenc}
\usepackage{lmodern,xcolor}
\usepackage{amsmath,amssymb,amsthm}
\usepackage{microtype}
\usepackage{enumitem}
\usepackage{authblk}
\usepackage{thmtools,thm-restate}
\usepackage[hidelinks]{hyperref}
\usepackage[nameinlink,noabbrev,capitalise]{cleveref}

\newtheorem{theorem}{Theorem}
\newtheorem{lemma}[theorem]{Lemma}
\newtheorem{proposition}[theorem]{Proposition}
\newtheorem{corollary}[theorem]{Corollary}
\theoremstyle{remark}
\newtheorem{remark}[theorem]{Remark}
\newtheorem{question}[theorem]{Question}

\crefname{question}{Question}{Questions}
\Crefname{question}{Question}{Questions}
\newcommand{\End}{\operatorname{End}}
\newcommand{\Aut}{\operatorname{Aut}}
\newcommand{\dist}{\operatorname{dist}}

\newcommand{\tr}{\operatorname{tr}}

\newcommand{\Spec}{\operatorname{Spec}}
\newcommand{\core}{\operatorname{Core}}
\newcommand{\Sh}{\operatorname{Sh}}
\newcommand{\Cay}{\operatorname{Cay}}

\newcommand{\F}{\mathbb F}
\newcommand{\R}{\mathbb R}
\newcommand{\one}{\mathbf1}
\newcommand{\f}{\mathbf f}

\newcommand{\blfootnote}[1]{%
  \begingroup
  \renewcommand{\thefootnote}{}%
  \footnote{#1}%
  \addtocounter{footnote}{-1}%
  \endgroup
}

\title{Cores and endomorphisms  of distance-regular graphs\protect\\
of small diameter}

\author[a]{Hong-Jun Ge}
\author[b,c,d]{Kolja Knauer}
\author[a]{Jack H. Koolen}

\affil[a]{School of Mathematical Sciences,
University of Science and Technology of China,
Hefei 230026, People's Republic of China}

\affil[b]{%
  School of Mathematical Sciences,
  Hebei Key Laboratory of Computational Mathematics and Applications,
  Hebei Normal University,
  Shijiazhuang 050024, P.~R.~China
}

\affil[c]{Departament de Matem\`atiques i Inform\`atica,
Universitat de Barcelona, Barcelona, Spain}

\affil[d]{Centre de Recerca Matem\`atica, Barcelona, Spain}

\date{}

\begin{document}
\maketitle

\blfootnote{E-mail addresses:
  \href{mailto:gehj567@gmail.com}
       {\texttt{gehj567@gmail.com}} (H.-J. Ge),
  \href{mailto:kolja.knauer@ub.edu}
       {\texttt{kolja.knauer@ub.edu}} (K. Knauer),
  \href{mailto:koolen@ustc.edu.cn}
       {\texttt{koolen@ustc.edu.cn}} (J. H. Koolen).
}

\vspace{-2em}

\begin{abstract}
We prove that every non-bipartite antipodal distance-regular graph of diameter three or four is a pseudocore, meaning that every endomorphism is either an automorphism or its image is a complete subgraph. 
The diameter three case settles a conjecture of Geertsema, Godsil and Guo.
On the other hand, we construct distance-regular graphs that are not core-complete in every even diameter at least four,
answering a question of Geertsema, Godsil and Guo in the negative.
Finally, we exhibit a non-bipartite antipodal distance-regular
graph of diameter five that is not a pseudocore.
\end{abstract}

\section{Introduction}

The \emph{core} of a finite graph $X$, denoted by $\core(X)$, is an (up to isomorphism unique) smallest induced subgraph to which $X$ admits a homomorphism. 
We call $X$ itself a \emph{core} if $\core(X)\cong X$,   \emph{core-complete} if it is a core or its core is
a complete graph, and a \emph{pseudocore} if the image of every non-bijective endomorphism is a complete subgraph.
Every pseudocore is therefore core-complete.

Symmetry can impose strong restrictions on graph cores.
Cameron and Kazanidis~\cite{CameronKazanidis2008} proved that
nonedge-transitive graphs, and hence rank-three graphs,
are core-complete.
Godsil and Royle~\cite{GodsilRoyle2011} established
core-completeness for distance-transitive graphs and showed
that every triangle-free non-bipartite distance-regular graph
is a core.
Roberson~\cite{Roberson19} proved that every primitive strongly
regular graph is a pseudocore, settling and strengthening
the conjecture of Cameron and Kazanidis.
Using vector colorings, Godsil et al.~\cite[Theorem~3.14]{GRRSV2019}
proved that every non-bipartite antipodal distance-regular graph
of diameter three with fiber size two  is a core.

For distance-regular graphs in general, Van Dam, Koolen and
Tanaka~\cite[Problem~45]{VDKT2016} asked which graphs are
core-complete, noting that no counterexamples were known
at the time.
More recently, Geertsema, Godsil and Guo~\cite{GGG} extended
Roberson's homomorphism-matrix method to distance-regular graphs.
They asked whether every distance-regular graph is
core-complete~\cite[Problem~1.1]{GGG}, and conjectured that
a non-bipartite antipodal distance-regular graph of diameter
three has no endomorphism whose image has diameter
two~\cite[Conjecture~1.2]{GGG}.

Our main theorem establishes the pseudocore property
in diameters three and four, 
where the diameter-three case settles a conjecture of
Geertsema, Godsil and Guo.

\begin{restatable}{theorem}{maintheorem}\label{thm:main}
Every non-bipartite antipodal distance-regular graph
of diameter three or four is a pseudocore.
\end{restatable}

The proof uses a criterion based on the cosine sequence for the smallest eigenvalue, see \cref{thm:three-value}. 
For the graphs in \cref{thm:main}, the criterion applies in diameter four and in diameter three when the cosines at distances one and three coincide. 
In the remaining diameter-three case, we show that the
graph is a core. 
We also show that the criterion can hold only in diameter at most four, with equality only for antipodal graphs, see \cref{thm:diameter}.

On the other hand, we give a negative answer to~\cite[Problem~1.1]{GGG}, i.e., we show that the distance-regular \emph{Doob graphs} $D(m,0)$
 are not core-complete if $m\ge 2$, see \cref{prop:doob}. 
This in particular shows that core-completeness is not determined
by the intersection array.

Finally, we consider the coset graph of the dual of the
ternary Golay code, a non-bipartite antipodal distance-regular
graph of diameter five.
We show that it is core-complete but not a pseudocore,
see \cref{prop:golay}.
Thus \cref{thm:main}
fails in diameter five.

\section{Preliminaries}\label{sec:prelim}

All graphs are finite and simple, and are connected
unless otherwise stated. 
An \emph{endomorphism} of a graph $X$ is a homomorphism
from $X$ to itself.
We denote the set of endomorphisms of $X$ by $\End(X)$.
We use $\one$ for the all-ones vector and $J$ for the
all-ones matrix, with dimensions determined by the context;
$J_s$ denotes the $s\times s$ all-ones matrix.
 For a homomorphism $\phi:X\to Y$, we regard its image as the subgraph of $Y$ induced by $\phi(V(X))$. 
A \emph{retraction} onto a subgraph $C$ is a
homomorphism $X\to C$ that fixes every vertex of $C$. 
In this case $C$ is called a \emph{retract} of $X$.
Every graph retracts onto a copy of its core.
Two graphs $X$ and $Y$ are \emph{homomorphically equivalent}
if there are homomorphisms $X\to Y$ and $Y\to X$.
This holds if and only if their cores are isomorphic,
see~\cite[Lemma~6.2.3]{GodsilRoyle2001}.

Let $X$ be a distance-regular graph of diameter $D$, valency $k$, and
intersection array
$ \{b_0,b_1,\ldots,b_{D-1};c_1,c_2,\ldots,c_D\}.$
We set $b_D=c_0=0$ and $a_i=k-b_i-c_i$. 
For $x\in V(X)$ and $0\le i\le D$, let
$ \Gamma_i(x):=\{y\in V(X):\dist(x,y)=i\}.$
Write $A_i$ for the
distance-$i$ matrix, $A:=A_1$, $k_i:=|\Gamma_i(x)|$,  and $n=|V(X)|$. 
There are polynomials
$p_i$ such that $A_i=p_i(A)$, and $k_i b_i=k_{i+1}c_{i+1}$ for $0\le i<D$.
We refer to~\cite{BCN1989,VDKT2016} for background on
distance-regular graphs. We recall the following monotonicity
property of their intersection numbers.
\begin{lemma}[{\cite[Proposition~4.1.6(i)]{BCN1989}}]
\label{lem:intersection-monotonicity}
Let $X$ be a distance-regular graph of diameter $D$ and valency $k$.
Then its intersection numbers satisfy
\[
  k=b_0>b_1\ge\cdots\ge b_{D-1}>b_D=0
  \qquad\text{and}\qquad
  1=c_1\le c_2\le\cdots\le c_D\le k.
\]
\end{lemma}

Let $\tau$ be the smallest eigenvalue of $X$, of multiplicity $m$,
and let $E$ be the corresponding spectral idempotent multiplied
by $n/m$. Thus
$E=\frac nm
\prod_{\substack{\lambda\in\Spec(X)\\\lambda\neq\tau}}
\frac{A-\lambda I}{\tau-\lambda}.$
Since
$\mathbb{R}[A]=\operatorname{span}\{I,A_1,\ldots,A_D\}$,
the entry $E(x,y)$ depends only on $\dist(x,y)$.
In particular, its diagonal entries are all equal to $1$.
Moreover, $E$ is a Gram matrix of unit vectors such that
\begin{equation}\label{eq:normalized-E}
	E^2=\frac nmE,\qquad E\one=0,\qquad
	\operatorname{rank}E=m,\qquad \tr E=n.
\end{equation}

Define $w_0,\ldots,w_D$ by
$E(x,y)=w_{\dist(x,y)}.$
The sequence $(w_0,\ldots,w_D)$ is the
\emph{smallest-eigenvalue cosine sequence}.
It satisfies
\begin{equation}\label{eq:cosine-recurrence}
	c_iw_{i-1}+a_iw_i+b_iw_{i+1}=\tau w_i
	\qquad(0\le i\le D),
\end{equation}
where the terms involving $w_{-1}$ and $w_{D+1}$ are omitted.
In particular,
\begin{equation}\label{eq:wi}
w_0=1,\qquad w_1=\frac{\tau}{k},\qquad
p_i(\tau)=k_iw_i,\qquad |w_i|\le1.
\end{equation}
For further details on these standard facts,
see~\cite[Sections~2.2 and~4.1]{BCN1989}.

Following Roberson~\cite[Definition~3.1]{Roberson19},
for $\phi\in\End(X)$, define the pullback matrix $M^\phi$
and its difference $\Delta^\phi$ from $E$ by
\[
 M^\phi(x,y):=w_{\dist(\phi(x),\phi(y))},
 \qquad \Delta^\phi:=M^\phi-E.
\]
The following lemma follows from \cite[Lemmas~3.1 and~3.3]{GGG} and the identities $A_i=p_i(A)$.

\begin{lemma}\label{lem:pullback}
For every $\phi\in\End(X)$ and $0\le i\le D$, the matrix $M^\phi$ is positive semidefinite and satisfies
\[
 M^\phi(A-\tau I)=0,\qquad
 M^\phi A_i=p_i(\tau)M^\phi,\qquad
 \Delta^\phi A_i=p_i(\tau)\Delta^\phi.
\]
In particular, if $\Delta^\phi(u,v)=0$ for some $u,v\in V(X)$,
then
$\sum_{x\in\Gamma_i(v)}\Delta^\phi(u,x)=0$ 
for $0\le i\le D$.
\end{lemma}

\begin{corollary}
\label{cor:extremal}
Suppose $w_i$ is a minimum or maximum of $w_0,\ldots,w_D$.
If $\dist(x,y)=i$ and $\phi\in\End(X)$, then
$ w_{\dist(\phi(x),\phi(y))}=w_i.$
\end{corollary}

\begin{proof}
The $(x,x)$ entry of $M^\phi A_i=p_i(\tau)M^\phi$ gives
$ \frac1{k_i}\sum_{y\in\Gamma_i(x)}
 w_{\dist(\phi(x),\phi(y))}=w_i.$
An average equals an extremal value only if all its terms equal that value.
\end{proof}

We recall the following property of endomorphisms, proved by
Ne\v{s}et\v{r}il~\cite{Nesetril1971}.

\begin{lemma}\label{lem:local-injectivity}
	An endomorphism of a finite connected graph that is injective on
	each neighborhood is an automorphism.
\end{lemma}

A distance-regular graph $X$ of diameter $D$ is \emph{antipodal} if the
relation $\dist(x,y)\in\{0,D\}$ is an equivalence relation. 
Its classes, called \emph{antipodal fibers}, have a common size $r\ge2$.
The \emph{antipodal quotient} $\overline X$ has the antipodal fibers as
vertices, with two distinct fibers adjacent whenever an edge of $X$
joins them. 
Let $\pi:V(X)\to V(\overline X)$ denote the quotient map. 
In this setting
$k_D=r-1$, $c_D=k$, and $a_D=0$.
We shall use the following standard covering property of
antipodal distance-regular graphs.

\begin{lemma}[{\cite[Theorem~2.1]{GodsilHensel92}}]
	\label{lem:antipodal-cover}
	Let $X$ be an antipodal distance-regular graph with diameter
	$D\ge3$, valency $k$, and antipodal fiber size $r$.
	Then $X$ is an $r$-fold cover of its antipodal quotient
	$\overline X$, meaning that the edges between any two distinct fibers
	are either empty or form a perfect matching.
	If $D=3$, then $\overline X\cong K_{k+1}$.
\end{lemma}

We also recall the following property of the smallest-eigenvalue
cosine sequence.

\begin{lemma}[{\cite[Corollary~4.1.2]{BCN1989}}]\label{lem:cosine-signs}
Let $X$ be a distance-regular graph of diameter $D$.  
Then $(-1)^i w_i>0$ for $0\le i\le D$.
\end{lemma}

\begin{lemma}[{\cite[Lemma~2.2]{Pascasio2001}}]\label{lem:even-cosines}
Let $X$ be a non-bipartite distance-regular graph  with diameter $D\ge 3$. Then $-1<w_i<1$ for $1\le i<D$.
Moreover, $w_D=1$ if and only if $X$ is antipodal and $D$ is even.
\end{lemma}
Consequently, if $X$ is a non-bipartite antipodal distance-regular
graph of even diameter $D\ge4$, then
$p_D(\tau)=k_Dw_D=r-1,$
and $E(x,y)=1$ if and only if $x$ and $y$ lie in the same
antipodal fiber.

\section{A three-value cosine criterion}\label{sec:three-value}

A strongly regular graph is called {\emph{primitive}} if both it and its complement
are connected.

\begin{lemma}[{\cite[Corollary 3.6]{Roberson19}}]\label{lem:prim srg}
	 Every primitive strongly regular graph is a pseudocore.
\end{lemma}

We now give a criterion that treats the diameter-three case with $w_1=w_3$ and the antipodal diameter-four case by the same argument.

\begin{theorem}
\label{thm:three-value}
Let $X$ be a distance-regular graph,
and let $(w_0,\ldots,w_D)$ be its smallest-eigenvalue
cosine sequence.
Suppose that 
$ \{w_0,w_1,\ldots,w_D\}=\{1,w_1,w_2\}$ and $ w_1<0<w_2<1.$
Then $X$ is a pseudocore.
\end{theorem}

\begin{proof}
	Let $\alpha:=w_1<0$ and $\beta:=w_2>0$.
	Define a relation on $V(X)$ by
	$
	x\equiv y$ if $E(x,y)=1.$
	Choose unit vectors $\mathbf{u}_x$ with $E(x,y)=\langle \mathbf{u}_x,\mathbf{u}_y\rangle$.
	Note that $E(x,y)=1$ holds exactly when $\mathbf{u}_x=\mathbf{u}_y$.
Thus $\equiv$ is an equivalence relation.
Every equivalence class has size
$ f=\sum_{i:\,w_i=1}k_i.$
Let $\widehat E$ be obtained from $E$ by retaining one row and
column from each class.
After ordering the vertices by equivalence classes, we have
$E=\widehat E\otimes J_f$.
Let $N=n/f$.
 By~\eqref{eq:normalized-E},
\begin{equation}\label{eq:quotient-Gram}
 \widehat E^2=\frac{N}{m}\widehat E,
 \qquad \widehat E\one=0,
 \qquad \operatorname{rank}\widehat E=m, \qquad \tr \widehat{E}=N.
\end{equation}
It follows that $\Spec(\widehat{E})=\{[0]^{N-m},[\frac{N}{m}]^m\}$.

Define a graph $Y$ on these classes by declaring two distinct
classes adjacent when the corresponding entry of $\widehat E$
is $\alpha$. 
Let $B$ be the adjacency matrix of $Y$.
Then
\begin{equation}\label{eq:two-inner-products}
 \widehat E=(1-\beta)I+(\alpha-\beta)B+\beta J.
\end{equation}
The identity $\widehat E\one=0$ implies that $Y$ is
$\frac{\beta N+1-\beta}{\beta-\alpha}$-regular.
Choosing an orthogonal eigenbasis of $\widehat E$ containing
$\one$ and applying Equation~\eqref{eq:two-inner-products}, we obtain
$$\Spec(Y)=\{[\frac{\beta N+1-\beta}{\beta-\alpha}],[\frac{1-\beta}{\beta-\alpha}]^{N-m-1},[\frac{N/m-1+\beta}{\alpha-\beta}]^m\}.$$
Since $w_2\neq1$, vertices at distance two in $X$ give distinct nonadjacent vertices of $Y$.
Thus $Y$ is not complete.
Consequently, $Y$ is a connected
strongly regular graph. 
It is primitive because a connected imprimitive
strongly regular graph is complete multipartite and has no
positive nonprincipal eigenvalue.

Let $\pi:X\to Y$ map vertices to their equivalence classes.
Let $\phi\in\End(X)$. 
Since $1$ and $\alpha$ are the extreme
cosine values, \cref{cor:extremal} shows that $\phi$ preserves
the relations $E(x,y)=1$ and $E(x,y)=\alpha$. It therefore induces
an endomorphism $\overline\phi$ of $Y$ satisfying
\begin{equation}\label{eq:pi chi}
	\pi\circ\phi=\overline\phi\circ\pi.
\end{equation}
By \cref{lem:prim srg}, either
$\overline\phi$ is an automorphism or its image is a complete subgraph of $Y$.

Suppose first that $\overline\phi$ is an automorphism. Distinct
vertices $u,v\in\Gamma_1(x)$ are at distance one or two, so
$E(u,v)\in\{\alpha,\beta\}$ and $\pi(u)\neq\pi(v)$.
Equation~\eqref{eq:pi chi} implies that $\phi(u)\neq\phi(v)$.
Thus $\phi$ is injective on every neighborhood.
\cref{lem:local-injectivity} implies that
$\phi$ is an automorphism.

Suppose next that the image of $\overline\phi$ is a clique $T$ of
$Y$. 
Since $w_2=\beta\notin\{1,\alpha\}$, no two vertices of
$\pi^{-1}(T)$ are at distance two in $X$.
 Hence every connected
component of $X[\pi^{-1}(T)]$ is a clique. The image of $\phi$
is connected and contained in this subgraph by Equation~\eqref{eq:pi chi}, so it is a clique.
\end{proof}

The three-value condition in \cref{thm:three-value} also
restricts the diameter of $X$, which may be of independent interest.

\begin{theorem}\label{thm:diameter}
Let $X$ be a distance-regular graph, and let
$(w_0,\ldots,w_D)$ be its smallest-eigenvalue cosine sequence.
Suppose that
$ \{w_0,w_1,\ldots,w_D\}=\{1,w_1,w_2\}
 $ and $w_1<0<w_2<1.$
Then $D\le4$, and if $D=4$, then $X$ is antipodal.
\end{theorem}

\begin{proof}
Let $\alpha=w_1$ and $\beta=w_2$. 
The graph $X$ is non-bipartite, since a
bipartite distance-regular graph has smallest-eigenvalue cosines
$w_i=(-1)^i$. 
We assume that $D\ge4$.
By \cref{lem:cosine-signs} and \cref{lem:even-cosines},
the smallest-eigenvalue cosines alternate in sign, and they satisfy
$-1<w_i<1$ for $1\le i<D$.
Consequently,
$ w_3=\alpha$ and $w_4\in\{\beta,1\}.$
We first show that $w_4=\beta$ is impossible. Suppose otherwise.
Comparing~\eqref{eq:cosine-recurrence} at $i=1$ and $i=3$ gives
\begin{equation}\label{eq:t}
 t:=\frac{1-\beta}{\beta-\alpha}=a_1-a_3>0.
\end{equation}

 Fix vertices $u,v$ with $\dist(u,v)=3$,
and partition $\Gamma_1(v)$ into
\[
 C=\Gamma_1(v)\cap\Gamma_2(u),\qquad
 S=\Gamma_1(v)\cap\Gamma_3(u),\qquad
 B=\Gamma_1(v)\cap\Gamma_4(u).
\]
Then $|C|=c_3$, $|S|=a_3$, and $|B|=b_3$.
Since $D\ge4$, both $C$ and $B$ are nonempty.
No edge joins $C$ to $B$.

Let $L$ be the adjacency matrix of the local graph $X[\Gamma_1(v)]$,
which is $a_1$-regular. Its vertices have mutual distance at most two
in $X$, so the principal submatrix of $E$ indexed by $\Gamma_1(v)$ is
\[
 E[\Gamma_1(v)]=(1-\beta)I+(\alpha-\beta)L+\beta J\succeq0,
\]
where $J$ is the all-ones matrix. Thus, for every
$\f\in\R^{\Gamma_1(v)}$ with $\f\perp\one$,
$\f^{\mathsf T}L\f\le t\|\f\|^2$. By~\eqref{eq:t},
\begin{equation}\label{eq:lower}
 \f^{\mathsf T}(a_1I-L)\f\ge a_3\|\f\|^2.
\end{equation}

For disjoint sets $U,W$, let $e(U,W)$ denote the number of edges
between them.
Now take $\f$ equal to $b_3$ on $C$, $-c_3$ on $B$, and $0$ on $S$.
Then $\f\perp\one$. 
Since there are no edges between $C$ and $B$,
the Laplacian quadratic form gives
\begin{align}
 \f^{\mathsf T}(a_1I-L)\f
 =b_3^2e(C,S)+c_3^2e(B,S)\le b_3^2c_3|S|+c_3^2b_3|S|
 =|S|\|\f\|^2,\label{eq:upper}
\end{align}
Comparing~\eqref{eq:lower} and~\eqref{eq:upper}, equality
must hold throughout. As $b_3,c_3>0$, we have $e(C,S)=c_3|S|$ and
$e(B,S)=b_3|S|$. 
Thus every vertex of $S$ is adjacent to every vertex of
$C\cup B$.

If $S$ is not empty, then the local degree gives
\[
 a_1\ge c_3+b_3=k-|S|,\qquad b_1=k-a_1-1\le |S|-1.
\]
On the other hand, any vertex of $C$ and any vertex of $B$ are
nonadjacent and have all of $S\cup\{v\}$ as common neighbours.
Therefore $c_2\ge |S|+1$, and hence $b_1<c_2$.
This contradicts $b_1\ge c_2$.
Indeed, any geodesic
$x_0,x_1,x_2,x_3$ satisfies
$ \Gamma_1(x_1)\cap\Gamma_1(x_3)
 \subseteq\Gamma_1(x_1)\cap\Gamma_2(x_0)$,
and these two sets have sizes $c_2$ and $b_1$, respectively.

It follows that $a_3=|S|=0$. The recurrence at $i=3$ now gives
$k\beta=\tau\alpha=k\alpha^2$, so $\beta=\alpha^2$.
Equation~\eqref{eq:t} yields
\(
 a_1=-\frac{1+\alpha}{\alpha}\ge1.
\)
At $i=2$ and $i=4$, the recurrence becomes
\[
 (k-a_i)\alpha+a_i\alpha^2=k\alpha^3.
\]
Here, for $i=4$, we use $b_4=0$ if $D=4$ and $w_5=\alpha$ otherwise.
Consequently,
\[
 a_2=a_4=k(1+\alpha)=\frac{ka_1}{a_1+1}\ge\frac{k}{2}.
\]
By \cref{lem:intersection-monotonicity}, we obtain
\[
 b_2=k-a_2-c_2<k-a_2\le a_2=a_4\le k-c_4\le k-c_3=b_3,
\]
a contradiction.

Therefore $w_4=1$. 
\Cref{lem:even-cosines} implies that $D=4$ and $X$
is antipodal. 
\end{proof}

\begin{corollary}\label{cor:D4}
Every antipodal distance-regular graph of diameter four is
core-complete. Every non-bipartite such graph is a pseudocore.
\end{corollary}

\begin{proof}
If $X$ is bipartite, then its core is $K_2$. Assume that $X$ is
non-bipartite. By \cref{lem:cosine-signs} and \cref{lem:even-cosines},
we have $0<w_2<1$ and $w_4=1$. Since $c_4=k$ and $a_4=0$,
the recurrence~\eqref{eq:cosine-recurrence} at $i=4$ gives
$kw_3=\tau w_4=\tau$. Hence $w_3=\tau/k=w_1<0$, and
\cref{thm:three-value} shows that $X$ is a pseudocore.
\end{proof}

\begin{corollary}\label{cor:fusion}
	Let $X$ be a distance-regular graph of diameter three whose
	smallest-eigenvalue cosines satisfy $w_1=w_3$.
	Then $X$ is core-complete.
	If $X$ is non-bipartite, then $X$ is a pseudocore.
\end{corollary}

\begin{proof}
If $X$ is bipartite, then its core is $K_2$.
Otherwise, \cref{lem:cosine-signs} and \cref{lem:even-cosines}
give $w_1=w_3<0<w_2<1$. 
Thus $X$ is a pseudocore by
\cref{thm:three-value}.
\end{proof}

\begin{remark}
The equality $w_1=w_3$ does occur: the line graph of the Petersen graph
is a non-bipartite antipodal distance-regular graph of diameter three with
this property. It is a counterexample to the cosine-injectivity assertion
in~\cite[Lemma~2.3]{GGG}, which is used in the proof
of~\cite[Lemma~4.4]{GGG}. 
Further non-bipartite antipodal distance-regular graphs of
diameter $3$ satisfying $w_1=w_3$ are provided by
Fon-der-Flaass~\cite[Construction~4]{FDF02}.
\Cref{cor:fusion} treats this case
without that assertion.
\end{remark}

\section{Antipodal graphs of diameter three}\label{sec:diameter-three}
In this section, let $X$ denote a non-bipartite antipodal distance-regular graph of diameter $3$ and antipodal fiber size $r\geq2$. 
By \cref{lem:antipodal-cover}, $X$ is an $r$-fold
cover of $K_{k+1}$.
Writing $c=c_2$, its intersection array is
$ \{k,(r-1)c,1;1,c,k\}.$
The characteristic polynomial of its intersection matrix factors as
$ (x-k)(x+1)\bigl(x^2+(rc-k+1)x-k\bigr).$
Thus its distinct eigenvalues are $k,\sigma,-1,\tau$, where $\sigma>0$
and $\tau$ is the smallest eigenvalue, with
\begin{equation}\label{eq:d3-spectrum}
 \sigma\tau=-k,\qquad \sigma+\tau=k-rc-1.
\end{equation}
As a  consequence,
\begin{equation}\label{eq:d3-a1}
	a_1=k-(r-1)c-1
	=\sigma-1+c\left(1-\frac r{\sigma+1}\right).
\end{equation}

Since $X$ is antipodal, $A_3=I_{n/r}\otimes (J_r-I_r).$
It follows that 
$\Spec(A_3)=\{[r-1]^{n/r},[-1]^{{n-n/r}}\}$.
By \eqref{eq:wi}, $w_1=\tau/k=-1/\sigma$ and $p_3(\tau)=k_3w_3$ is an eigenvalue of $A_3$ and $k_3=r-1$.
By \cref{lem:cosine-signs}, we have $w_3<0$,
so $p_3(\tau)=(r-1)w_3=-1$.
 Using $k_2=k(r-1)$ and $E\one=0$
gives
\begin{equation}\label{eq:d3-cosines}
 w_0=1,\qquad w_1=-\frac1\sigma,\qquad
 w_2=\frac1{(r-1)\sigma}>0,\qquad w_3=-\frac1{r-1}.
\end{equation}
As $p_3(\tau)=-1$, by \cref{lem:pullback}, every $\phi\in\End(X)$ satisfies
\begin{equation}\label{eq:d3-fiber-sum}
 M^\phi(A_3+I)=0.
\end{equation}

We first show that an endomorphism preserving one antipodal pair also preserves the corresponding two rows of the smallest-eigenvalue matrix.
The following argument adapts the sign argument in the proof of
\cite[Theorem~3.4]{Roberson19}.

\begin{lemma}\label{lem:delta}
Let $\phi\in\End(X)$, and suppose that
$ \dist(u,v)=\dist(\phi(u),\phi(v))=3.$
Then the $u$- and $v$-rows of $\Delta^\phi=M^\phi-E$ are zero.
\end{lemma}

\begin{proof}
Let $F$ be the fiber containing $u$ and $v$.
If $x\sim v$,
then
$\dist(u,x)=\dist(\phi(u),\phi(x))=2$,
since $D=3$ and graph homomorphisms do not increase distances.
It follows that $\Delta^\phi(u,x)=0$ for every
$x\in\Gamma_1(v)$.

Since $\Delta^\phi(u,v)=0$, by \cref{lem:pullback},
$
 \sum_{x\in\Gamma_2(v)}\Delta^\phi(u,x)=0.$
 The terms with $x\sim u$ vanish because $\phi$ preserves adjacency.
Every other $x$ in
this sum has distance $2$ from both $u$ and $v$. Moreover,
$\phi(x)\neq\phi(u)$, since $\dist(\phi(v),\phi(x))\leq2$ whereas
$\dist(\phi(v),\phi(u))=3$. 
Its contribution is therefore one of
$\{0, w_1-w_2, w_3-w_2.\}$
The last two quantities are strictly negative by~\eqref{eq:d3-cosines}.
Thus $\Delta^\phi(u,x)=0$ for every
$x\in\Gamma_2(v)$.

For $x\in F$, all neighbours of $x$ lie outside $F$. Hence, by \cref{lem:pullback},
$\Delta^\phi A=\tau\Delta^\phi$ gives
\[
 \tau\Delta^\phi(u,x)=\sum_{y\sim x}\Delta^\phi(u,y)=\sum_{\substack{y\sim x\\y\in \Gamma_1(v)\cup \Gamma_2(v)}}\Delta^\phi(u,y)=0.
\]
As $\tau\neq0$, the $u$-row is zero. Interchanging $u$ and $v$ proves the
same assertion for the $v$-row.
\end{proof}

\begin{theorem}\label{thm:diam3-rigidity}
If the smallest-eigenvalue cosines of $X$ satisfy $w_1\neq w_3$, then $X$ is a core.
\end{theorem}

\begin{proof}
Let $\phi\in\End(X)$.
Suppose that there is a vertex $x$ in a fiber $F$ such that
$\dist(\phi(x),\phi(y))\le 2$ for every $y\in F$. 
For $i=0,1,2$, let
$ n_i=\bigl|\{y\in F:\dist(\phi(x),\phi(y))=i\}\bigr|.$
Then $ n_0+n_1+n_2=r$.
Since $F=\{x\}\cup\Gamma_3(x)$, by Equation~\eqref{eq:d3-fiber-sum}, taking the $(x,x)$-entry of
$M^\phi(A_3+I)=0$ gives
\[
0=\sum_{y\in F}M^\phi(x,y)
=\sum_{y\in F}w_{\dist(\phi(x),\phi(y))}
=n_0w_0+n_1w_1+n_2w_2.
\]
Equation~\eqref{eq:d3-cosines} gives
$  \sigma n_0=n_1-\frac{n_2}{r-1}.$
Since $n_0\ge1$, the eigenvalue $\sigma$ is rational and hence an integer.
Eliminating $n_1$ yields
$ n_2=\frac{(r-1)\bigl(r-n_0(\sigma+1)\bigr)}r.$
Its nonnegativity gives $0<n_0(\sigma+1)\leq r$, while its integrality
and $\gcd(r,r-1)=1$ give $r\mid n_0(\sigma+1)$. 
Therefore
$ r=n_0(\sigma+1)$ and $ n_2=0.$
Since $n_2=0$, we have
$\phi(F)\subseteq\{\phi(x)\}\cup\Gamma_1(\phi(x)).$
Thus any two vertices of $\phi(F)$ are at distance at most
two in $X$.
Repeating the argument for every vertex of $F$ shows that
$\phi(F)$ is a clique of size $\sigma+1$, with constant
preimage multiplicity $s=r/(\sigma+1)$.
An edge of this clique has at least $\sigma-1$ common
neighbours. By~\eqref{eq:d3-a1},
$ \sigma-1\leq a_1=\sigma-1-c(s-1).$
Thus $s=1$, giving $\sigma=r-1$ and hence $w_1=w_3$,
contrary to the hypothesis.

Therefore, for each vertex $x$ there is a vertex $y$ such that $\dist(\phi(x),\phi(y))=3$.
It follows that $\Delta^\phi=0$ by \cref{lem:delta}.
The cosines $w_0,w_1,w_2,w_3$ are pairwise distinct.
Thus $\phi$ is an automorphism.
\end{proof}

\begin{remark}
When $r=2$, the equality $w_1=w_3$ would give $w_1=-1$
and hence $\tau=-k$, forcing $X$ to be bipartite.
Thus \cref{thm:diam3-rigidity} recovers the result that
every non-bipartite Taylor graph is a core
\cite[Theorem~3.14]{GRRSV2019}.    
\end{remark}

\begin{corollary}\label{cor:diam3}
Every antipodal distance-regular graph of diameter $3$ is core-complete.
Every non-bipartite such graph is a pseudocore.
\end{corollary}

\begin{proof}
	If $X$ is bipartite, then its core is $K_2$.
Otherwise, if $w_1=w_3$, apply
\cref{cor:fusion}; if $w_1\neq w_3$, apply
\cref{thm:diam3-rigidity}.
\end{proof}

\maintheorem*

\begin{proof}
This follows from \cref{cor:D4,cor:diam3}.
\end{proof}

\section{Examples and further questions}\label{sec:examples}

\subsection{Pseudocores}
Consider the following
 matrix over $\F_3$:
\[
H=\left[\begin{array}{cccccc|ccccc}
1&0&0&0&0&0&1&1&1&1&1\\
0&1&0&0&0&0&1&1&-1&-1&0\\
0&0&1&0&0&0&1&-1&1&0&-1\\
0&0&0&1&0&0&-1&1&0&1&-1\\
0&0&0&0&1&0&-1&0&1&-1&1\\
0&0&0&0&0&1&0&-1&-1&1&1
\end{array}\right].
\]

Writing $h_1,\ldots,h_{11}$ for the columns of $H$, we define
$\Gamma:=\operatorname{Cay}(\mathbb{F}_3^6,
\{\pm h_i:1\leq i\leq 11\})$.
The row space of $H$ is an $[11,6,5]_3$ code and hence, up to
equivalence, the \emph{ternary Golay code}~\cite{Pless1992}.
Thus $H$ is a parity-check matrix of its dual, and $\Gamma$
is isomorphic to the coset graph of the dual ternary Golay code.
It is a non-bipartite antipodal distance-regular graph of
diameter five with intersection array
$\{22,20,18,2,1;1,2,9,20,22\}$.
Blokhuis, Brouwer and Haemers~\cite[Section~3.4]{BBH2007}
showed that $\Gamma$ is $3$-chromatic and uniquely determined
by its intersection array.
Moreover, $\Gamma$ is a distance-transitive $3$-fold antipodal
cover of the Berlekamp--van Lint--Seidel
graph~\cite[Proposition~3.31]{Alfuraidan2013}.
We show below that $\Gamma$ is not a pseudocore, hence, \cref{thm:main} fails in diameter five.

\begin{proposition}\label{prop:golay}
The graph $\Gamma$ 
is not a pseudocore, but core-complete.
\end{proposition}

\begin{proof}
Let $e_1,\ldots,e_6$ be the standard basis of $\F_3^6$, put
$U=\operatorname{Span} \{e_1,e_2\}$, and define the linear map
\[
 \rho(x_1,\ldots,x_6)
 =(x_1+x_3+x_4+x_5,\ x_2+x_6,\ 0,0,0,0).
\]
Note that for $1\le i\le 11$, $\rho(h_i)\in\{e_1,e_2\}=\{h_1,h_2\}$.
Thus $\rho$  is an endomorphism of $\Gamma$. 
It fixes $U$ pointwise, and
$\{\pm h_i\}\cap U=\{\pm e_1,\pm e_2\}$. 
It follows that $\rho$ is a retraction onto
\[
 \Gamma[U]\cong
 \Cay(\F_3^2,\{\pm e_1,\pm e_2\})\cong H(2,3).
\]
So $\Gamma$ is not a pseudocore. Since $\Gamma$ is $3$-chromatic~\cite[Section~3.4]{BBH2007} and contains a triangle,
$\operatorname{Core}(\Gamma)\cong K_3$. Hence, $\Gamma$ is core-complete.
\end{proof}

We wonder whether \cref{thm:main} holds beyond diameter five.

\begin{question}\label{q:antipodal-pseudocore}
Is every non-bipartite antipodal distance-regular graph
of diameter at least six a pseudocore?
\end{question}

\subsection{Core-completeness}
We first describe the cores of the Doob graphs.
This yields distance-regular graphs that are not core-complete
in every even diameter at least four, addressing the existence
issue noted in~\cite[Problem~45]{VDKT2016} and giving a negative
answer to~\cite[Problem~1.1]{GGG}.

The \emph{Shrikhande graph} $\Sh$ is the Cayley graph
$\Cay(\mathbb Z_4^2,S)$, where
$S=\{\pm(1,0),\pm(0,1),\pm(1,1)\}$.
For integers $m,n\ge0$ with $m+n\ge1$, the \emph{Doob graph}
$D(m,n)$ is the Cartesian product
\[
 D(m,n):=\Sh^{\square m}\square K_4^{\square n}.
\]
These graphs are distance-regular of diameter $2m+n$, with the same
intersection arrays as the Hamming graphs $H(2m+n,4)$;
see~\cite[Section~9.2B]{BCN1989}.
Since $\core(H(2m+n,4))\cong K_4$, 
the following proposition shows that core-completeness is not determined by the intersection array.

\begin{proposition}\label{prop:doob}
The Doob graphs have cores
\[
 \core(D(m,n))\cong
 \begin{cases}
 \Sh,&n=0,\ m\ge1,\\
 K_4,&n\ge1.
 \end{cases}
\]
Consequently, $D(m,n)$ is not core-complete if and only if $n=0$ and $m>1$.
\end{proposition}

\begin{proof}
The map
$(x_1,\ldots,x_m)\longmapsto
 (x_1+\cdots+x_m,0,\ldots,0)$
is a retraction of $D(m,0)$ onto a copy of $\Sh$.
Since $\Sh$ is a core, see \cite{MPRR}, this proves the first case.

If $n\ge1$, varying a single $K_4$ coordinate gives a copy of $K_4$ in $D(m,n)$. 
The chromatic number of a Cartesian product is the maximum of the chromatic numbers of its factors. 
Since $\chi(\Sh)=\chi(K_4)=4$, it follows that $\chi(D(m,n))=4$.
Thus $D(m,n)$ and $K_4$ are homomorphically equivalent, proving the second case. The final assertion follows immediately.
\end{proof}

Taking $m\ge2$ and $n=0$ gives a graph of every even diameter
$2m\ge4$ whose core is the noncomplete graph $\Sh$.
This leads to the following question.

\begin{question}\label{q:odd-diameter}
Is every distance-regular graph of odd diameter core-complete?
\end{question}

Although the diameter-five graph $\Gamma$ in \cref{prop:golay}
is not a pseudocore, it is core-complete.
Together with \cref{cor:D4,cor:diam3}, this example motivates the following question.

\begin{question}\label{q:antipodal-core}
Is every antipodal distance-regular graph core-complete?
\end{question}

\subsection{Endomorphism monoids}
Under composition, $\End(X)$ is a monoid whose group
of units is the automorphism group $\Aut(X)$.
Let $M=\End(X)$ and let $C$ be a copy of the core of $X$.
By standard results on finite transformation
monoids~\cite[Section~3.1]{Steinberg2010}, the minimal ideal
$I$ of $M$ consists of the endomorphisms of minimum image
size and is completely simple. If $\rho$ is a retraction
onto $C$, then $\rho\in I$, and restriction to $C$ yields
$\rho M\rho\cong\operatorname{Aut}(C)$. Hence, by the
Rees--Suschkewitsch theorem~\cite[Theorem~3.3.1]{Howie1995},
$I$ is a Rees matrix semigroup over $\operatorname{Aut}(C)$.

Thus, if $X$ is core-complete, then $M$ is either a group
or its minimal ideal is a Rees matrix semigroup over
the symmetric group $S_{|V(C)|}$.
If $X$ is a pseudocore but not a core,
then every noninvertible endomorphism belongs to $I$.
Consequently, $M$ is completely regular (a union of groups)
with exactly two $\mathcal{J}$-classes, namely
$\operatorname{Aut}(X)$ and $I$.
A finite monoid is called \emph{aperiodic} if all its
subgroups are trivial.
In particular, the only aperiodic monoid arising as the
endomorphism monoid of a core-complete graph with an edge
is the trivial monoid.

Knauer and Puig i Surroca~\cite{KnauerPuig2025}
proved that every finite monoid is isomorphic to the
endomorphism monoid of a regular graph.
Since distance-regular graphs need not be core-complete,
we repeat the question raised in
\cite[Section 6]{KnauerPuig2025}.

\begin{question}
For every finite monoid $M$, does there exist a
distance-regular graph $X$ such that
$\operatorname{End}(X)\cong M$?
\end{question}

\section*{Acknowledgements}

We thank Hebei Normal University, and especially Sergey Goryainov,
for their hospitality and support during our visit.

Kolja Knauer is partially supported by grant PID2022-137283NB-C22,
funded by MICIU/\allowbreak AEI/\allowbreak
10.13039/\allowbreak 501100011033 and ERDF/EU;
by the Severo Ochoa and Mar\'ia de Maeztu Program for Centers
and Units of Excellence in R\&D (CEX2020-001084-M);
and by grant ANR-21-CE48-0012  and the grant of The Natural Science Foundation of Hebei Province (project No. A2023205045).

Jack H.~Koolen is supported by the National Natural Science
Foundation of China (No.~12471335) and the Anhui Initiative
in Quantum Information Technologies (No.~AHY150000).

\paragraph{Statement of AI use.}
The authors obtained the diameter-three result and the Doob graph
counterexamples. During subsequent work on the diameter-four case,
ChatGPT 5.6 Sol suggested that core-completeness lifts from the
antipodal quotient to a non-bipartite antipodal distance-regular
graph of even diameter. The authors developed this observation
into a lifting argument proving the diameter-four case.
They then recognized that this argument shares the method used
in the diameter-three case with $w_1=w_3$, leading them to formulate
and prove the more general \cref{thm:three-value}.
Since the diameter-four result follows directly from this theorem,
the separate lifting theorem was removed from the manuscript.
The authors take full responsibility for the mathematical claims, proofs, and final presentation of the paper.

\small
\bibliographystyle{my-siam}
\bibliography{antipodal_core_complete}

\end{document}